# An Efficient Hybrid CS and K-Means Algorithm for the Capacitated P-Median Problem

H G Mazinan [a], G R Ahmadi [b], E Khaji [c]

[a] *School of Railway Engineering, Iran University of Science and Technology, 1684613114 Narmak,Tehran, Iran;*
[b] *Department of Industrial Engineering, Amir Kabir University, Tafresh, Iran;* [c] *Master Student in Complex Adaptive Systems, Department of Physics, University of Göteborg, 41296 Göteborg, Sweden.*

Capacitated p-median problem (CPMP) is an important variation of facility location problem in which p capacitated medians are economically selected to serve a set of demand vertices so that the total assigned demand to each of the candidate medians must not exceed its capacity. This paper surveys and analyses the combination of Cuckoo Search and K-Means algorithms to solve the CPMP. In order to check for quality and validity of the suggestive method, we compared the final solution produced over the two test problems of Osman and Christofides, each of which including 10 sample tests. According to the results, the suggested meta-heuristic algorithm shows superiority over the rest known algorithms in this field as all the best known solutions in the first problem set, and several sample sets in the second problem set have been improved within reasonable periods of time.

___________________________________________________________________________

## 1. Introduction

The capacitated P-median problem (CPMP) is an NP-complete problem which investigates the problem of partitioning a set of N nodes into M different disjoint clusters, each represented by a certain node that is designed as concentrator. The N–M nodes that are not chosen as concentrators are referred as terminals. The partitioning of the initial N nodes must be performed in such a way that a measure of total distance between the terminals and their corresponding concentrators is minimized. In addition, a capacity constraint imposed on the concentrators must be met, in order to obtain feasible solutions to the problem [1–4]. A direct application of the CPMP is in the context of communication networks deployment, where a set of terminals in the network must be assigned to the corresponding concentrator in order to compose access networks that satisfy the rate requirements of such terminals [5]. In this context, most of the efforts so far has focused on the topological design of communication networks (e.g. Wireless Sensor Networks (WSN), backbone networks or mobile networks [6–8]) since many of the processes involved in such networks can be approached as a CPMP problem, e.g. server load balancing in backbone links, and clustering and/or data aggregation in WSN. This problem also stems from applications in other fields, such as the design of distribution networks [9], vehicle routing [10], political district organization [11], selection of facilities for a university's admission examination [2] and sales force territories design [12], among others. To efficiently tackle this NP-complete problem, a large amount of algorithms can be found in the literature, e.g. see [1, 2, 13–15, 3, 16–19] and references therein. Among such previous approaches, evolutionary-based algorithms have been intensively investigated [2, 5, 20, 21]. These algorithms usually perform well in discrete and extremely constrained optimization problems such as the here considered CPMP, and consequently they are generally shown to outperform other algorithmic approaches in terms of complexity. In spite of this huge previous work on evolutionary computing for the CPMP, the performance of grouping-based heuristic algorithms has been showed to be efficient. In this context, the most widely used grouping approach in the literature is the grouping genetic algorithm (GGA), which essentially is a class of evolutionary algorithms specially tailored for grouping problems, i.e. problems where a number of items must be assigned to a set of predefined groups. These series of algorithms have also been discussed in many literatures [21-34, 39].

This paper presents a novel hybrid meta-heuristic approach which is a combination of CS and K-means algorithm. As the results will clearly show, the proposed approach is, to the authors' knowledge, the best evolutionary techniques for the CPMP problem reported so far in the literature.

The remainder of the paper is structured as follows: Section 2 poses the CPMP from a mathematical standpoint. Next, Section 3 comprises the technical content of this contribution where the proposed hybrid algorithm is comprehensively described. Special emphasis is placed on key features of such algorithm, and brief descriptions about its components such as K-means algorithm, cuckoo search, and Roulette wheel algorithm. Section 4 presents the experimental part of the paper, where the performance of the proposed algorithm is assessed in CPMP instances of diverse size. Section 5 concludes the paper by drawing some final remarks.

## 2. The mathematical model

The CPMP considered in this paper is modeled as a binary integer-programming problem (P):

$$v(p) = Min \sum_{i \in V} \sum_{j \in V} d_{ij} x_{ij}, \tag{1}$$
$$S.t: \sum_{j \in V} x_{ij} = 1, \quad \forall i \in V, \tag{2}$$
$$x_{ij} \leq y_{ij}, \quad \forall i, j \in V, \tag{3}$$
$$\sum_{j \in V} y_j = p, \tag{4}$$
$$\sum_{i \in N} q_i x_{ij} \leq Q_j y_j, \quad \forall j \in V, \tag{5}$$
$$x_{ij}, y_j \in \{0,1\}, \quad \forall i \in V, j \in V, \tag{6}$$

where $V = \{1,2,\ldots,n\}$ is a set of vertices and possible medians in the graph, $d_{ij}$ are distance values, $x_{ij}$ is allocation decision variable with $x_{ij}=1$, if demand point i is allocated to median j, and $x_{ij}=0$ otherwise. $y_j$ is the median decision variable with $y_j=1$, if a facility is located at point j, and $y_j=0$, otherwise. Finally $q_i$ and $Q_j$ are demand and the capacity of points i and j respectively.

The objective function (1) is to minimize the sum of distances between demand points and medians. Constraint (2) imposes each customer to be assigned to exactly one median. Disaggregated inequality (3) guarantees that each demand point must be allocated to a point that has been selected as median. Constraint (4) ensures that always p medians are selected, and constraint (5) imposes that the sum of demands assigned to a median does not exceed its capacity.

## 3. The proposed hybrid CS & K-means technique

This work is aiming at utilizing a combination of cuckoo search and k-means in solving CPMP. In the following sections we discus briefly the so called algorithms, and moreover, how the proposed algorithm have combined these approaches in order to obtain a homogenous technique without the disadvantageous

### 3.1. The traditional K-means algorithm

K-means is one of the most familiar learning algorithms that solve the well known clustering problem. The procedure follows a simple way to classify a given data set through a certain number of clusters of k fixed a priori. The main problem is to define k-centroids, one for each cluster. These centroids should be placed in a cunning way since different location causes different result. Therefore, the first step of this algorithm is to determine k random points as the centroids of clusters, preferably, far away from each other. Thence, each point belonging to a given data set is assigned to the nearest centroid. When no point is pending, the first step is completed and an early groupage is done. At this point, it needs to re-calculate k new centroids as barycenters of the clusters resulting from the previous step. Having generated these k new centroids, a new binding has to be done between the same data set points and the nearest new centroid which constituted a convergent loop in the algorithm. Technically speaking, the algorithm aims at minimizing an objective function, in this case a squared error function. The objective function

$$J = \sum_{j=1}^{k} \sum_{i=1}^{n} \left| x_i^{(j)} - c_j \right|^2, \qquad (7)$$

where $\left|x_i^{(j)} - c_j\right|^2$ is a chosen distance measure between a data point and the cluster center $c_j$, is an indicator of the distance of the n data points from their respective cluster centers.

Although it can be proved that the procedure will always terminate, the k-means algorithm does not necessarily find the most optimal configuration, corresponding to the global objective function minimum. The algorithm is also significantly sensitive to the initial randomly selected cluster centers. The k-means algorithm can be run multiple times to reduce this effect.

## 3.2. Cuckoo Search

In this section, the breeding behavior of cuckoos and the characteristics of Levy flights in reaching best habitat societies are overviewed. Cuckoo Search algorithm (CS) is one of the well known heuristic algorithms [35, 36] by Xin-She Yang and Suash Deb in 2009. It has been developed by simulating the intelligent breeding behavior of cuckoos. It is a population-based search procedure that is used as an optimization tool, in solving complex, nonlinear and non-convex optimization problems [37]. Female parasitic cuckoos specialize and lay eggs that closely resemble the eggs of their chosen host nest. Cuckoo chooses this host nest by natural selection. The shell of the cuckoo egg is usually thick. They have two different layers with an outer chalky layer which resists to cracking when the eggs are dropped in the host nest. The cuckoo egg hatches earlier than the host bird's egg, and the cuckoo chick grows faster. Thence, remaining cuckoo's egg in the nest are recognized by host birds with a probability $P_a \in [0, 1]$ and these eggs are thrown away or the nest is abandoned, and a completely new nest is built, in a new location by host bird. The mature cuckoo form societies and each society have its habitat region to live in. The best habitat from all of the societies will be the destination for the cuckoos in other societies. Thence, they immigrate toward this best habitat. A randomly distributed initial population of host nest is generated and then the population of solutions is subjected to repeated cycles of the search process of the cuckoo birds to lay an egg. The cuckoo randomly chooses the host nest position to lay an egg using Levy flights random-walk and is given in Equations. (8) and (9).

$$V_{pq}^{t+1} = V_{pq}^t + S_{pq} \times Levy(\lambda) \times \alpha, \tag{8}$$

$$Levy(\lambda) = \frac{\Gamma(1 + \lambda) \times \sin\left(\pi \times \frac{\lambda}{2}\right)}{\Gamma\left(\frac{1+\lambda}{2}\right) \times \lambda \times 2^{\frac{\lambda-1}{2}}}^{1/\lambda}, \tag{9}$$

where λ is a constant (1 ≤ λ ≤ 3), t is the current generation number and α is a random number generated between [−1, 1]. Also s > 0 is the step size which should be related to the scales of the problem of interests. If s is a big number, then, the new generated solution will have a long distance with the old solution. Therefore, such a move is unlikely to be accepted. If it is too small, the change is too small to be significant, and consequently such search is not efficient. So, a proper step size is important to maintain the search as efficiently as possible [37]. Hence, the step size is calculated using Eq. (22).

$$S_{pq} = V_{pq}^t - V_{fq}^t, \qquad (10)$$

where $p, f \in \{1, 2, \ldots, m\}$, and $q \in \{1, 2 \ldots D\}$ are randomly chosen indexes, and f and p are not equal. D is the number of parameters to be optimized and m is the total population of host nest positions. Utilizing (20), the cuckoo chooses the host nest or communal nest, and an egg laid by a cuckoo is evaluated. The host bird identifies the alien egg with the probability value associated with the egg's fitness value using Eq. (23).

$$Pro_p = (0.9 * Fit_p / max(Fit)) + 0.1, \qquad (11)$$

where $Fit_p$ is the fitness value of the solution p which is proportional to the fitness value of an egg in the nest position p, and $Pro_p$ gives the survival probability rate of the cuckoo's egg. If the random generated probability $p_a \in [0, 1]$ is greater than $Pro_p$, then the alien egg is identified by the host bird. Then, the host bird destroys the alien egg away, or abandons the nest, and cuckoo finds a new position using Equation (24) for laying an egg. Otherwise, the egg grows up and will be alive for the next generation based on the fitness function.

$$x_p = x_{p\,min} + rand\,(0, 1) \times (x_{p\,max} - x_{p\,min}), \qquad (12)$$

where $x_{p\,min}$ and $x_{p\,max}$ are the minimum and maximum limits of the parameter to be optimized.

### 3.3. Roulette Wheels

Fitness proportionate selection, also known as roulette-wheel selection, is a genetic operator used in genetic algorithms for selecting potentially useful solutions for recombination.

In fitness proportionate selection, as in all selection methods, the fitness function assigns a value to each possible solution or chromosomes. This value is used to associate a probability of selection with each individual chromosome. If $f_i$ is the fitness of individual $i$ in the population, its probability of being selected is

$$p_i = \frac{f_i}{\sum_{j=1}^{N} f_i}. \tag{13}$$

N is the number of individuals in the population. Usually a proportion of the wheel is assigned to each of the possible selections based on their fitness value. This could be achieved by dividing the fitness of a selection by the total fitness of all the selections, thereby normalizing them to 1. Then a random selection is made similar to how the roulette wheel is rotated.

While candidate solutions with a higher fitness will be less likely to be eliminated, there is still a chance that they may be. Contrast this with a less sophisticated selection algorithm, such as truncation selection, which will eliminate a fixed percentage of the weakest candidates. With fitness proportionate selection there is a chance some weaker solutions may survive the selection process; this is an advantage, as though a solution may be weak, it may include some component which could prove useful following the recombination process. In the proposed algorithm, we have utilized this advantageous to increase the possibility of converging to global minimum in K-means algorithm.

### 3.4. Generation of the population

Generating an initial population with perfect fitness can ease the process of finding the global minimum within reasonable time in any meta-heuristic hybrid algorithm. In this paper, we have utilized a method to generate well-fitted individuals (nests) which always guarantees to have a good enough population within all the generations. This approach includes three different steps, *Local improvement* (K-means Clustering), optimization of generation (CS), and *mutation*.

Local improvement is done using K-means clustering. Clustering the nodes, determining a median in each cluster, and assigning the nodes of each cluster to their clusters' median are included in this step. However, the assigning part should be passed according to the following routine; in each cluster, after choosing a median among all the nodes, the nodes are assigned to that median according to the distance between each node and the median. The node which has the lowest distance to the median is in the first priority, and the remaining nodes have to be assigned to the median in their cluster according to their distance. Finally, all the distances are calculated for all the clusters, and one can find a value for each cluster which determines the quality or the fitness of that cluster. It should be considered that by the end of this step, we have not generated any individual for the population yet, and it should be finalized in the second step.

In the next part, the median of each cluster must change its location, and again, all the nodes in that cluster must be assigned to the new median. This step continues until the best fitness value would be obtained for all the clusters. It should be noted that the fitness value of the total clusters does depend on the node that is first assigned to the median in its cluster. Depending on the median from which the assigning procedure starts, a different fitness value is obtained. Furthermore, any infeasible generated nest (individual) in which the summation of the nodes demands in a cluster (clusters) is more than their median capacity should be removed in the very moment after its generation. Having recognized all the clusters in a (nest), each median in a cluster should change until the best possible fitness value for those clusters is obtained. In conclusion, a nest (a set of clusters among several possible cluster sets) is locally improved.

Obviously, by changing the order of the node in which we start the k-means algorithm we can get the other possible cluster sets (nest). If all of these nests are locally improved, thence an initial population with a quit good quality is prepared. It should be noticed that each new individual which is generated not only in the initial population, but in all the generations should be locally improved to enter the algorithm.

### 3.5. Mutation

In order to avoid getting stuck in local minima, in each generation, the mutation procedure must be followed. This step includes choosing k nodes among those which have the worst fitness within all clusters, and reassigning them to the other medians according to the K-means and Roulette Wheel algorithms. Thence, the new individuals should be again locally improved regarding the so called procedure in the previous step. As it was mentioned in section Roulette Wheel, the authors took advantage of the fact that using Roulette Wheel gives opportunity to even the worst alternatives to have a chance to be selected in a random competition. Considering the fact that K-means algorithm has the disadvantageous of getting stuck in local minima in non convex problems such as CPMP, using Roulette Wheel in selection section of K-means helps to perturb the algorithm, and increase the number of examined answers in the feasible set. A scheme of the proposed hybrid algorithm can be seen in Figure 1.

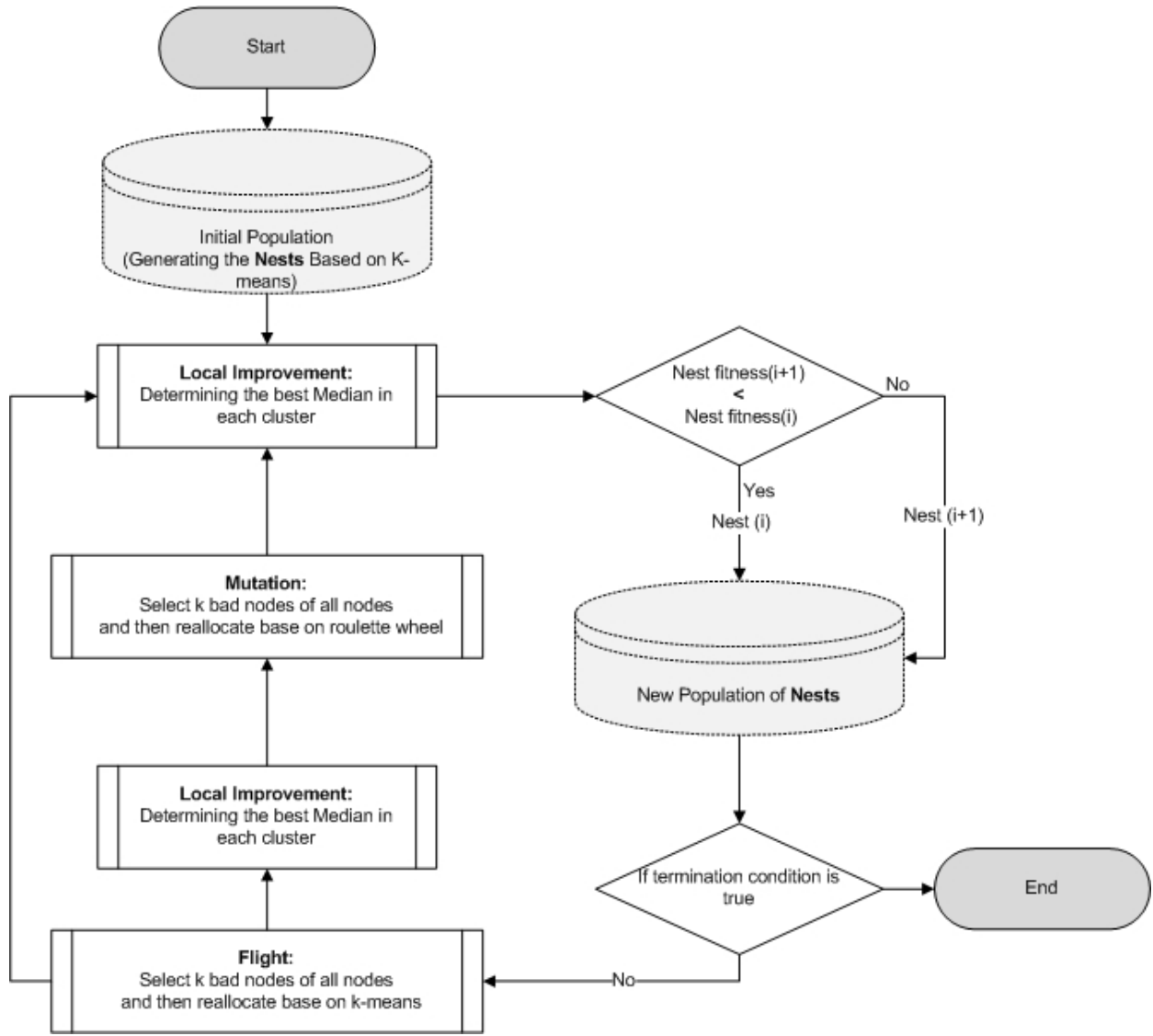

Figure 1: The scheme of CS & K-means algorithm.

## 3.6. New Population of nests

New population of nests should be generated in each generation according to the so called procedure in CS algorithm. The algorithm is briefly discussed in section 3.2.

## 4. Computational experiments and results

In order to evaluate the performance of the proposed technique when applied to the CPMP, we have carried out experiments and comparisons with existing evolutionary-based approaches by utilizing examples of Osman and Christofides [38]. The samples include 50 nodes in the network, and the capacity of each median is equal to 120. The numbers of medians for all samples in the first and second problem sets are considered 5 and 10 respectively. Also, the demand amount of each vertex is different in each sample.

Specifically, we compare the proposed hybrid algorithm with the former proposed algorithm by the authors in [40]: the first includes CS and K-means while the second is a combination of GA and Ant Colony optimization. Each of the so called algorithms utilizes a meta-heuristic algorithm to optimize each individual (nest) in each generation (CS and GA respectively), and the next algorithm to assign nodes to the medians (K-means and Ant Colony algorithm).

The problem sets for the simulations can be obtained from http://people.brunel.ac.uk/mastjjb/jeb/info.html. It should be considered that the best computational results for these 20 test problems are calculated using Euclidean distances and reported in [38] after they are rounded down to the nearest integer. Computational experiments were carried out on a 3.2 GHz processor and 2 GB RAM. The main code of the proposed algorithm is implemented by JAVA. The initial values for running of GACO are the same as the ones in [40]. The simulation results which are the best and average solutions after 500 generations over the first and second problem sets of Osman and Christofides are reported in Table 1 and Table 2. Moreover, the scheme of the best result of simulation over sample 10 of the first problem set, which is achieved by CS and K-means algorithm, is shown in Figure 2.

**Table 1:** Comparison of the results obtained by the proposed hybrid CS and K-means algorithm in the first set of experiments run in Osman and Christofides examples.

| No. | medians | Best-Known | CS & K-means | | | | | GACO | | |
|---|---|---|---|---|---|---|---|---|---|---|
| | | | Best fitness | Time (s) | Dev. | | | Best fitness | Time (s) | % Dev. |
| | | | | | Improvement | Gap | | | | |
| 1 | 5 | 713 | 690 | 0.9 | 23.3 | 3.4% | - | - | 713 | 6.2 | | |
| 2 | 5 | 740 | 730 | 8.8 | 9.5 | 1.3% | - | - | 740 | 5.4 | | |
| 3 | 5 | 751 | 731 | 6.8 | 19.7 | 2.7% | - | - | 751 | 7.0 | | |
| 4 | 5 | 651 | 637 | 1.8 | 13.8 | 2.2% | - | - | 651 | 9.5 | | |
| 5 | 5 | 664 | 663 | 2.0 | 0.9 | 0.1% | - | - | 664 | 15.0 | | |
| 6 | 5 | 778 | 765 | 2.7 | 13.3 | 1.7% | - | - | 778 | 7.8 | | |
| 7 | 5 | 787 | 758 | 2.8 | 28.9 | 3.8% | - | - | 787 | 11.0 | | |
| 8 | 5 | 820 | 784 | 2.2 | 35.6 | 4.5% | - | - | 820 | 10.1 | | |
| 9 | 5 | 715 | 691 | 2.0 | 23.9 | 3.5% | - | - | 715 | 14.2 | | |
| 10 | 5 | 829 | 763 | 2.0 | 66.2 | 8.7% | - | - | 829 | 17.6 | | |

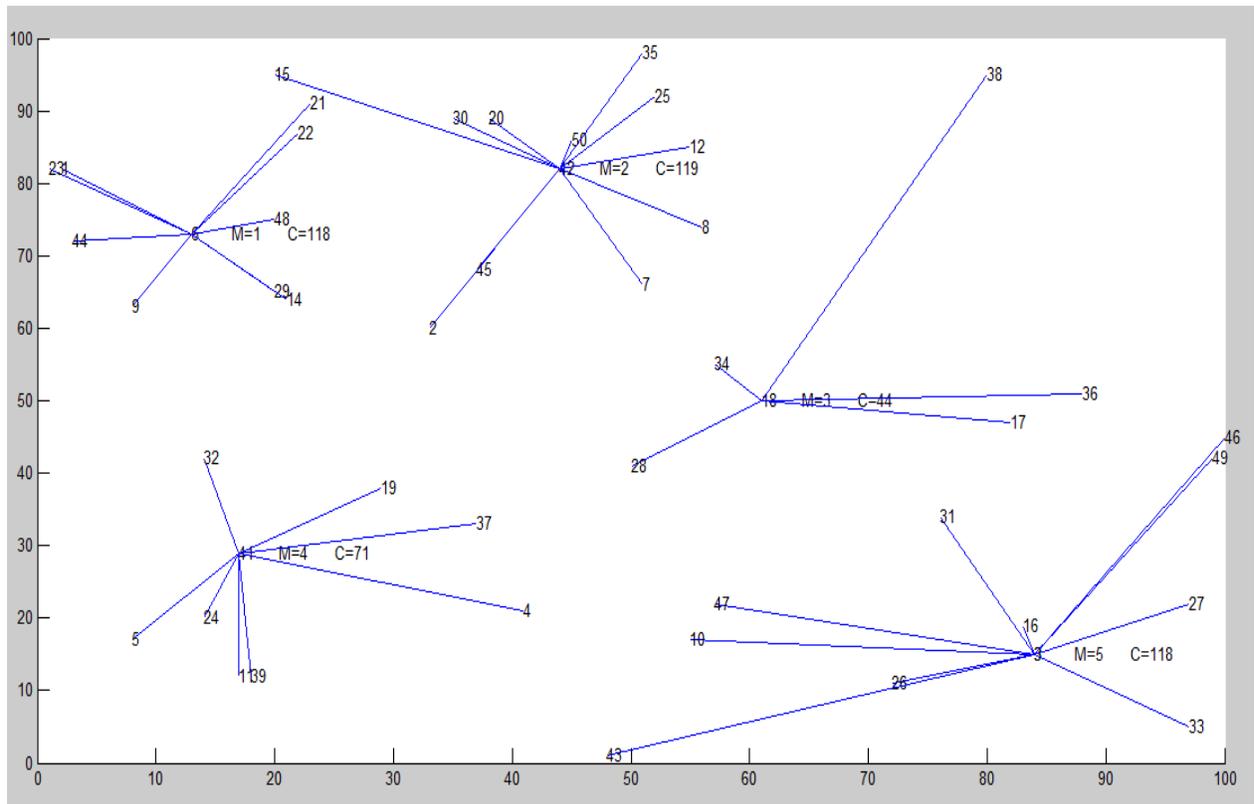

Figure 2: The best solution for sample 10 in the first problem set of Osman and Christofides (number of nodes=50, number of medians=5, best achieved fitness by CS and K-means=763).

Table 2: Comparison of the results obtained by the proposed hybrid CS and K-means algorithm in the second set of experiments run in Osman and Christofides examples.

| No. | medians | Best-Known | CS & K-means | | | | | GACO | | | |
|---|---|---|---|---|---|---|---|---|---|---|---|
| | | | Best fitness | Time (s) | Dev. | | | Best fitness | Time (s) | | % Dev. |
| | | | | | Improvement | | Gap | | | | |
| 1 | 10 | 1006 | 1006 | 49.6 | - | - | - | - | 1006 | 25.4 | | |
| 2 | 10 | 966 | 954 | 10.3 | 12.3 | 1.3% | - | - | 968 | 34.0 | 2 | 0.2% |
| 3 | 10 | 1026 | 1029 | 47.0 | - | - | 3 | 0.3% | 1026 | 29.4 | | |
| 4 | 10 | 982 | 984 | 51.7 | - | - | 2 | 0.3% | 986 | 38.3 | 4 | 0.4% |
| 5 | 10 | 1091 | 1082 | 15.2 | 9.3 | 0.9% | - | - | 1096 | 94.7 | 5 | 0.5% |
| 6 | 10 | 954 | 938 | 4.2 | 16.2 | 1.7% | - | - | 957 | 74.9 | 3 | 0.3% |
| 7 | 10 | 1034 | 1035 | 50.3 | - | - | 1 | 0.1% | 1034 | 98.1 | | |
| 8 | 10 | 1043 | 1019 | 8.0 | 24.0 | 2.4% | - | - | 1050 | 72.9 | 7 | 0.7% |
| 9 | 10 | 1031 | 1030 | 81.4 | 1.4 | 0.1% | - | - | 1037 | 89.1 | 6 | 0.6% |
| 10 | 10 | 1005 | 988 | 3.2 | 17.5 | 1.8% | - | - | 1022 | 140.6 | 17 | 1.7% |

According to Table 1, the proposed method has produced perfect solutions with superiority over all best known results including the former algorithm of the authors GACO. The degree of improvement varies within 0.9 to 66 (0.1% to 8.69%) in all samples of the first problem set. Regarding to Table 2, perfect solutions have been produced which are superior over many of the best known results. The degree of improvement varies within 1.4 to 24 (0.1% to 2.4%) in 6 out of 10 sample sets whereas the gap value differs from 1 to 3 (0.1% to 0.3%) in 4 data sets. Therefore, Comparison of results signifies that the proposed work is thoroughly superb as it is able to improve several best known results in the so-called problem sets in quite reasonable periods of time. Additionally, detailed information about the experiments over the first and second problem sets such as number of assigned nodes to each median in the best achieve results and number of demand units each median must supply are presented in Tables 3 and 4.

Table 3: Detail Information about the results of simulations over the problem set 1 achieved by CS & K-means.

| No. of Sample | No of assigned nodes to each median | No of demands which each of medians should supply | No. of Sample | No of assigned nodes to each median | No of demands which each of medians should supply |
|---|---|---|---|---|---|
| **1** | | | **6** | | |
| Median1 | 12 | 89 | Median1 | 11 | 111 |
| Median2 | 5 | 46 | Median2 | 9 | 105 |
| Median3 | 11 | 88 | Median3 | 9 | 97 |
| Median4 | 9 | 91 | Median4 | 11 | 105 |
| Median5 | 13 | 115 | Median5 | 10 | 90 |
| **2** | | | **7** | | |
| Median1 | 9 | 89 | Median1 | 10 | 112 |
| Median2 | 11 | 65 | Median2 | 11 | 111 |
| Median3 | 10 | 102 | Median3 | 12 | 94 |
| Median4 | 11 | 111 | Median4 | 8 | 81 |
| Median5 | 9 | 88 | Median5 | 9 | 91 |
| **3** | | | **8** | | |
| Median1 | 6 | 57 | Median1 | 10 | 113 |
| Median2 | 7 | 80 | Median2 | 12 | 120 |
| Median3 | 14 | 116 | Median3 | 10 | 85 |
| Median4 | 12 | 110 | Median4 | 7 | 68 |
| Median5 | 11 | 69 | Median5 | 11 | 106 |

| | 4 | | | 9 | |
|---|---|---|---|---|---|
| Median1 | 12 | 105 | Median1 | 9 | 70 |
| Median2 | 7 | 41 | Median2 | 8 | 87 |
| Median3 | 14 | 117 | Median3 | 12 | 120 |
| Median4 | 8 | 80 | Median4 | 11 | 117 |
| Median5 | 9 | 95 | Median5 | 10 | 91 |
| | 5 | | | 10 | |
| Median1 | 9 | 109 | Median1 | 10 | 118 |
| Median2 | 14 | 104 | Median2 | 12 | 119 |
| Median3 | 9 | 96 | Median3 | 6 | 44 |
| Median4 | 8 | 105 | Median4 | 9 | 71 |
| Median5 | 10 | 82 | Median5 | 13 | 118 |

Table 4: Detail Information about the results of simulations over the problem set 2 achieved by CS & K-means.

| No. of Sample | No of assigned nodes to each median | No of demands which each of medians should supply | No. of Sample | No of assigned nodes to each median | No of demands which each of medians should supply |
|---|---|---|---|---|---|
| | 1 | | | 6 | |
| Median1 | 13 | 104 | Median1 | 13 | 118 |
| Median2 | 10 | 81 | Median2 | 13 | 103 |
| Median3 | 11 | 116 | Median3 | 14 | 110 |
| Median4 | 9 | 90 | Median4 | 7 | 89 |
| Median5 | 8 | 59 | Median5 | 12 | 119 |
| Median6 | 6 | 86 | Median6 | 6 | 44 |
| Median7 | 13 | 119 | Median7 | 6 | 53 |
| Median8 | 10 | 113 | Median8 | 12 | 108 |
| Median9 | 8 | 72 | Median9 | 7 | 85 |
| Median10 | 12 | 89 | Median10 | 10 | 108 |
| | 2 | | | 7 | |
| Median1 | 9 | 86 | Median1 | 10 | 91 |
| Median2 | 8 | 60 | Median2 | 8 | 87 |
| Median3 | 8 | 59 | Median3 | 7 | 74 |
| Median4 | 11 | 89 | Median4 | 9 | 112 |
| Median5 | 6 | 52 | Median5 | 12 | 107 |
| Median6 | 13 | 110 | Median6 | 15 | 120 |
| Median7 | 8 | 98 | Median7 | 10 | 58 |
| Median8 | 13 | 115 | Median8 | 11 | 118 |
| Median9 | 11 | 120 | Median9 | 12 | 112 |
| **Median10** | 13 | 109 | **Median10** | 6 | 55 |
| **3** | | | **8** | | |
| **Median1** | 9 | 71 | **Median1** | 11 | 111 |
| **Median2** | 13 | 116 | **Median2** | 8 | 67 |
| **Median3** | 12 | 110 | **Median3** | 10 | 116 |
| **Median4** | 12 | 104 | **Median4** | 8 | 59 |
| **Median5** | 10 | 71 | **Median5** | 12 | 118 |
| **Median6** | 10 | 115 | **Median6** | 12 | 115 |
| **Median7** | 9 | 94 | **Median7** | 9 | 77 |
| **Median8** | 8 | 69 | **Median8** | 9 | 93 |
| **Median9** | 7 | 80 | **Median9** | 10 | 104 |
| **Median10** | 10 | 94 | **Median10** | 11 | 106 |
| **4** | | | **9** | | |
| **Median1** | 14 | 105 | **Median1** | 9 | 92 |
| **Median2** | 8 | 68 | **Median2** | 8 | 61 |
| **Median3** | 10 | 112 | **Median3** | 10 | 114 |
| **Median4** | 12 | 101 | **Median4** | 12 | 109 |
| **Median5** | 9 | 106 | **Median5** | 11 | 95 |
| **Median6** | 6 | 32 | **Median6** | 9 | 100 |
| **Median7** | 11 | 112 | **Median7** | 13 | 120 |
| **Median8** | 8 | 95 | **Median8** | 9 | 100 |
| **Median9** | 13 | 111 | **Median9** | 11 | 107 |
| **Median10** | 9 | 108 | **Median10** | 8 | 53 |
| **5** | | | **10** | | |
| **Median1** | 9 | 94 | **Median1** | 9 | 102 |
| **Median2** | 10 | 101 | **Median2** | 13 | 117 |
| **Median3** | 10 | 110 | **Median3** | 12 | 113 |
| **Median4** | 12 | 82 | **Median4** | 7 | 71 |
| **Median5** | 9 | 63 | **Median5** | 11 | 113 |
| **Median6** | 9 | 94 | **Median6** | 11 | 119 |
| **Median7** | 13 | 92 | **Median7** | 12 | 105 |
| **Median8** | 10 | 92 | **Median8** | 8 | 69 |
| **Median9** | 9 | 54 | **Median9** | 9 | 117 |
| **Median10** | 9 | 118 | **Median10** | 8 | 76 |

As the results indicates, the numbers of assigned nodes to each median are clearly near each other except the first and last examples whereas the numbers of demands which each of the medians should supply have more considerable differences. This signifies that in order to reach to the best global answer, the demand value of a client is a more significant factor compared with its location in determination of the median to which the node should be assigned.

# 5. Conclusions

In this paper we have studied the performance of hybrid CS & K-means algorithm when applied to the capacitated P-median problem (CPMP). The algorithm starts with generating initial population using K-Means algorithm. Thereafter, the new generation is resulted by applying the Cuckoo Search CS algorithm followed by Levy Flight. Thence, each new nest should be *Local Improved* and *Mutated* according to K-Means and Roulette Wheel algorithms. The new generation is generated according to CS. The algorithm continues until the termination condition would be satisfied. We have assessed the performance of the proposed technique by addressing CPMP instances of very diverse size, and by comparing the obtained results with those of existing best known outcomes and also the former work of the authors. In all the performed experiments, the proposed method is shown to outperform significantly such previous evolutionary techniques. CS- K-means is, therefore, proven to be an advantageous option to tackle this problem in future applications.